\documentclass [a4paper] {book}
\usepackage [T1] {fontenc}
\usepackage [utf8] {inputenc}
\usepackage [english]{babel}
\usepackage{amsmath,amssymb}
\usepackage{tikz-cd}
\newcommand{\numberset}[1]{\mathbb{#1}}
\newcommand{\N}{\numberset{N}}
\newcommand{\Z}{\numberset{Z}}

\newcommand{\C}{\numberset{C}} 
 \usepackage[a4paper,top=2cm,bottom=2cm,left=2.5cm,right=2.5cm,%
            heightrounded]{geometry}

\DeclareMathOperator{\Ker}{Ker}

\DeclareMathOperator{\reg}{reg}
\DeclareMathOperator{\sing}{sing}

\DeclareMathOperator{\AVF}{AVF}
\DeclareMathOperator{\IVF}{IVF}
\DeclareMathOperator{\Lie}{Lie}
\DeclareMathOperator{\alg}{alg}
\DeclareMathOperator{\SL}{SL}
\DeclareMathOperator{\hol}{hol}
\DeclareMathOperator{\VF}{VF}
\DeclareMathOperator{\PSL}{PSL}
\DeclareMathOperator{\GL}{GL}
\begin{document}
\large
\begin{center}
{\bfseries NEW EXAMPLES OF STEIN MANIFOLDS WITH VOLUME DENSITY PROPERTY}
\end{center} 
\normalsize
\begin{center}
Giorgio De Vito 
\end{center}
\begin{center}
Universit\"{a}t Bern \\ Mathematisches Institut (MAI) \\ Alpeneggstra\ss e 22 \\ 3012 Bern \\ Schweiz \\  \end{center}
\begin{center}
giorgio.devito@math.unibe.ch
\end{center}
{\scshape Abstract.} In the present paper we shall provide new examples of Stein manifolds enjoying the (algebraic) volume density property and compute their homology groups.
\large
\begin{center}
{\scshape 1. Introduction}
\end{center}
In 1997 J.--P. Rosay offered a highly valuable account of what is widely known as {\itshape Anders\'{e}n--Lempert theory} in his survey [R]. The main feature of this theory, based on Anders\'{e}n and Lempert's groundbreaking work (see [A] and [AL]), is that local injective holomorphic maps on holomorphically convex compact subsets can be approximated by global holomorphic automorphisms. Such a property was formalised via the concepts of ({\itshape algebraic}) {\itshape density property} and ({\itshape algebraic}) {\itshape volume density property}, introduced by Varolin ([V1]) and defined, according to [KKPS], as follows:\\ \\
{\bfseries Definition 1.} A complex manifold $X$ has the {\itshape density property} if, in the compact-open topology, the Lie algebra $\Lie_{\hol}(X)$ generated by completely integrable holomorphic vector fields on $X$ is dense in the Lie algebra $\VF_{\hol}(X)$ of all holomorphic vector fields on $X$. An affine algebraic manifold $X$ has the {\itshape algebraic density property} if the Lie algebra $\Lie_{\alg}(X)$ generated by completely integrable algebraic vector fields on it coincides with the Lie algebra $\AVF(X)$ of all algebraic vector fields on it.\\
On the other hand, suppose a complex manifold $X$ is equipped with a holomorphic volume form $\omega.$ We say that $X$ has the {\itshape volume density property with respect to $\omega$} if, in the compact-open topology, the Lie algebra $\Lie_{\hol}^{\omega}(X)$ generated by completely integrable holomorphic vector fields $\nu$ such that $\nu(\omega)=0$ is dense in the Lie algebra $\VF_{\hol}^{\omega}(X)$ of all holomorphic vector fields annihilating $\omega$. If $X$ is affine algebraic, we say that $X$ has the {\itshape algebraic volume density property with respect to an algebraic volume form $\omega$} if the Lie algebra $\Lie_{\alg}^{\omega}(X)$ generated by completely integrable algebraic vector fields $\nu$ such that $\nu(\omega)=0$ coincides with the Lie algebra $\AVF_{\omega}(X)$ of all algebraic vector fields annihilating $\omega.$\\ \\
We remark that algebraic density property and algebraic volume density property imply density property and volume density property respectively. For a proof of the second implication we refer to [KKAVDP]. \\
As we shall see, various examples of such highly symmetric objects were found over time and the above theory has been applied to many natural geometric questions, for an exhaustive list of which we refer to par. 4 of [K]. \\
\begin{center}
1.1. {\bfseries Notations and results}
\end{center} 
In this paper we will use the following notations, equivalent, as we shall explain, to those in [IK]. For $n\geq 3$ consider the family $\{M_n\}$ of $2\times 2$ matrices defined inductively as follows: \\

\noindent $M_3=$\;\( 
\begin{pmatrix}
1 & 0\\
z_1 & 1
\end{pmatrix}
\)
\( 
\begin{pmatrix}
1 & z_2\\
0 & 1
\end{pmatrix}
\)
\( 
\begin{pmatrix}
1 & 0\\
z_3 & 1
\end{pmatrix}
\)
=
 \( 
\begin{pmatrix}
1+z_2z_3 & z_2\\
z_1+z_3(z_1z_2+1) & z_1z_2+1
\end{pmatrix}
\),
\\ \\
$M_4=M_3$
\( 
\begin{pmatrix}
1 & z_4\\
0 & 1
\end{pmatrix}
\) =
\( 
\begin{pmatrix}
1+z_2z_3 & z_4(1+z_2z_3)+z_2\\
z_1+z_3(z_1z_2+1) & z_4(z_1+z_3(z_1z_2+1))+z_1z_2+1
\end{pmatrix}
\),
\\ \\
\dots \\ \\
$M_n=M_{n-1}$
\(
\begin{pmatrix}
1 & 0\\
z_n & 1
\end{pmatrix}
\)
\mbox{ if $n$ is odd}, \\ \\
$M_n=M_{n-1}$
\(
\begin{pmatrix}
1 & z_n\\
0 & 1
\end{pmatrix}
\)
\mbox{ if $n$ is even}.\\ \\
Also, for $n\geq 3$ consider the family $\{X_n\}$ of smooth complex algebraic hypersurfaces of $\C^n$ defined inductively as follows: 
$$X_3=\{(z_1,z_2,z_3)\in\C^3:(M_3)_{2,1}=1\}=\{p_3(z_1,z_2,z_3)=z_1+z_3+z_1z_3z_2-1=0\},$$ 
$$X_4=\{(z_1,\dots,z_4)\in\C^4:(M_4)_{2,2}=2\}=\{p_4(z_1,\dots,z_4)=z_1z_2-1+z_4(z_1+z_3+z_1z_3z_2)=0\},$$ $\dots$
$$X_n=\{(M_n)_{2,1}=1\}=\{p_n(z_1,\dots,z_n)=p_{n-2}(z_1,\dots,z_{n-2})+z_n(p_{n-1}(z_1,\dots,z_{n-1})+2)=0\}$$ if $n\geq 5$ is odd,
$$X_n=\{(M_n)_{2,2}=2\}=\{p_n(z_1,\dots,z_n)=p_{n-2}(z_1,\dots,z_{n-2})+z_n(p_{n-1}(z_1,\dots,z_{n-1})+1)=0\}$$ if $n\geq 6$ is even. \\ 
Furthermore, let us denote by $\Theta:\AVF_{\omega}(X_n)\to\mathcal{Z}_{n-2}(X_n)$ the isomorphism sending the $\omega$-divergence-free algebraic vector field $\xi$ to the closed $(n-2)$-form $\iota_{\xi}\omega$, defined via the interior product $\iota$ of $\xi$ and $\omega$, by $\Lie_{\alg}^{\omega}(X_n)$ the Lie algebra generated by the set $\IVF_{\omega}(X_n)$ of complete $\omega$-divergence-free algebraic vector fields on $X_n$ and by $\mathcal{B}_{n-2}(X_n)$ the space of exact algebraic $(n-2)$-forms on $X_n$. \\ \\
The results of this article follow. Proofs are given in sections 3.2 and 3.3. We just remark, at this stage, that the approach developed in [KKAV] has proven itself highly fruitful for the present work. \\ \\
{\bfseries Theorem 1.1.1.} The following statements hold true:\\ \\
1a) $X_3\subset\C^3_{z_1,z_2,z_3}$ is the affine modification, via the projection $\sigma:X_3\to\C^2_{z_1,z_3}$, of $\C^2$ along the divisor $\Delta=\{z_1z_3=0\}\subset\C^2$ with centre at $C_3=\{z_1z_3=1-z_1-z_3=0\}=\{(1,0),(0,1)\}\subset\reg\Delta\subset\Delta$.\\ \\
1b) $X_4\subset\C^4_{z_1,z_2,z_3,z_4}$ is the affine modification, via the projection $\sigma:X_4\to\C^3_{z_1,z_2,z_3}$, of $\C^3$ along the smooth divisor $X_3^0=\{z_1+z_3+z_1z_3z_2=0\}\subset\C^3$ with centre at $C_4=\{z_1z_2-1=z_1+z_3+z_1z_2z_3=0\}=\{(z_1,1/z_1,-z_1/2):z_1\in\C^*_{z_1}\} \simeq\C^*,\; C_4\subset X_3^0$. \\ \\
1c) If $n\geq 5$ is odd $X_n\subset\C^n_{z_1,\dots,z_n}$ is the affine modification, via the projection $\sigma:X_n\to\C^{n-1}_{z_1,\dots,z_{n-1}}$, of $\C^{n-1}$ along the smooth divisor $X_{n-1}^0=\{p_{n-1}(z_1,\dots,z_{n-1})+2=0\}\subset\C^{n-1}$ with centre at $C_n=\{p_{n-2}(z_1,\dots,z_{n-2})=p_{n-1}(z_1,\dots,z_{n-1})+2=0\}\subset X_{n-1}^0$, where $X_{n-1}^0\simeq\C^{n-2}\setminus X_{n-2}^0$ and $C_n\simeq X_{n-2}$.\\ 
If $n\geq 6$ is even $X_n\subset\C^n_{z_1,\dots,z_n}$ is the affine modification, via the projection $\sigma:X_n\to\C^{n-1}_{z_1,\dots,z_{n-1}}$, of $\C^{n-1}$ along the smooth divisor $X_{n-1}^0=\{p_{n-1}(z_1,\dots,z_{n-1})+1=0\}\subset\C^{n-1}$ with centre at $C_n=\{p_{n-2}(z_1,\dots,z_{n-2})=p_{n-1}(z_1,\dots,z_{n-1})+1=0\}\subset X_{n-1}^0$, where $X_{n-1}^0\simeq\C^{n-2}\setminus X_{n-2}^0$ and $C_n\simeq X_{n-2}$.\\ \\
2a) $H_0(X_3)=\Z,\;\;H_1(X_3)=0,\;\;H_2(X_3)=\Z ,\;\;H_k(X_3)=0\mbox{ for } k\geq 3.$\\ \\
2b) $H_0(X_4)=\Z,\;\;H_1(X_4)\cong H_2(X_4)=0,\;\;H_3(X_4)=\Z,\;\;H_k(X_4)=0\mbox{ for } k\geq4.$ \\ \\
2c) If $n\geq 5$ is odd, for $0\leq j\leq n-1$ the homology groups are alternately $\Z$ or $0$, starting from $H_0(X_n)=\Z$. In particular $H_{n-2}(X_n)=0$. \\ 
If $n\geq 6$ is even, for $0\leq j\leq \frac{n-2}{2}$ the homology groups are alternately $\Z$ or $0$ starting from $H_0(X_n)=\Z$, for $\frac{n}{2}\leq j\leq n-1$ the homology groups are alternately $\Z$ or $0$ starting from $H_{n-1}(X_n)=\Z$. In particular $H_{n-2}(X_n)=0$. \\ \\
Also, here comes our final \\ \\
{\bfseries Theorem 1.1.2.} For each $n\geq 3$ the variety $X_n\subset\C^n$ is such that $\Theta(\Lie_{\alg}^{\omega}(X_n))\supset\mathcal{B}_{n-2}(X_n)$, i.e. for every algebraic $(n-3)$-form $\alpha$ on $X_n$ there exists $\xi\in\Lie_{\alg}^{\omega}(X_n)$ such that $\Theta(\xi)=\iota_{\xi}\omega={\rm d}\alpha$. Hence $\Lie_{\alg}^{\omega}(X_n)=\AVF_{\omega}(X_n)$, i.e. $X_n$ has the algebraic volume density property. \\
\begin{center}
1.2. {\bfseries Motivations and purposes of our study}
\end{center} 
Relatively few manifolds are known to enjoy the volume density property so far. Apart from $\C^n$, the example $M^2:=\C^2\setminus\{xy=1\}$ equipped with $\omega=\frac{1}{xy-1}\;dx\wedge dy$ is due to Varolin (see [V]), along with $(\C^*)^n$ and $\SL_2(\C).$ \\
Kaliman and Kutzschebauch proved the algebraic volume density property in [KKAVDP] for $\SL_2(\C)$, for $\PSL_2(\C)$, for all semi-simple groups and for all linear algebraic groups with respect to the invariant volume form. Further, they extended this result to every connected affine homogenous space of a linear algebraic group $G$ over $\C$ which admits a $G$-invariant volume form (see [KKH]).
They also showed that the hypersurface $X':=\{P(u,v,\bar{x})=uv-p(\bar{x})=0\}\subset\C^{n+2}_{u,v,\bar{x}}$ (which is called {\itshape Danielewski surface} if $n=1$) enjoys, under some technical assumptions on the polynomial $p$, the algebraic volume density property with respect to an $\omega'$ s.t. $\omega'\wedge {\rm d}P=\Omega|_{X'},$ where $\Omega$ is the standard volume form on $\C^{n+2}$ (see [KKAVDP]). The holomorphic version of this result is due to Ramos-Peon (see [P]). Besides, as in [F], let us give the following \\ \\ 
{\bfseries Definition 1.2.} Let $Y$ be a Stein manifold.
\begin{enumerate}
\item $Y$ is {\itshape universal for proper holomorphic embeddings} if every Stein manifold $X$ with
$2\dim X < \dim Y$ admits a proper holomorphic embedding $X\hookrightarrow Y.$
\item $Y$ is {\itshape strongly universal for proper holomorphic embeddings} if, under the assumptions in (1), every continuous map $f_0 : X\to Y$ which is holomorphic in a neighborhood of a compact $\mathcal{O}(X)$-convex set $K\subset X$ is homotopic to a proper holomorphic embedding $f_0:X\hookrightarrow Y$ by a homotopy $f_t:X\to Y\;(t\in[0,1])$ such that $f_t$ is holomorphic and arbitrarily close to $f_0$ on $K$ for every $t\in[0, 1].$
\item $Y$ is (strongly) {\itshape universal for proper holomorphic immersions} if condition (1) (resp. (2)) holds for proper holomorphic immersions $X\to Y$ from any Stein manifold $X$ satisfying $2 \dim X\leq \dim Y.$ 
\end{enumerate}
It is well known that every Stein manifold with the density property is an Oka manifold and is Gromov-elliptic (see [FH]). The following theorem has been proven by Andrist and Wold (see [AW]), if $X$ is an open Riemann surface, by Andrist et al. (see [AF], Theorems 1.1-1.2) for embeddings and by Forstneri\u{c} (see [F2], Theorem 1.1) for immersions in the double dimension. \\ \\
{\bfseries Theorem 1.2.1.} Let $X$ be a Stein manifold with the density or the volume density property. Such $X$ is strongly universal for proper holomorphic embeddings and immersions. \\ \\
The above makes Stein manifolds with volume density property ideal target spaces for a positive solution of \\ \\ 
{\bfseries L\'{a}russon's question.} Does every Stein manifold admit an {\itshape acyclic} (i.e. being a homotopy equivalence) {\itshape embedding} into an Oka manifold? \\ \\
In other words, we are looking for a generalisation of the Remmert-Bishop-Narasimhan embedding theorem (see [N]). \\ 
We note T. Ritter has recently proven that every open Riemann surface acyclically embeds into an elliptic manifold. In case the surface is an annulus, he showed the wished target can be $\C\times\C^*$ (see [T], [T1]).\\ \\
The above observations encourage and motivate the present paper, in which we shall provide new examples of Stein manifolds enjoying the (algebraic) volume density property and compute their homology groups. \\ \\ We remark that our manifolds arise as fibres in a fibration Ivarsson and Kutzschebauch employ in [IK] to affirmatively solve the following \\ \\
{\bfseries Gromov--Vaserstein problem} (see [G]). Let $f:\C^n\to\SL_m(\C)$ be a holomorphic map. Does $f$ factorise as a finite product of holomorphic maps sending $\C^n$ into unipotent subgroups in $\SL_m(\C)$? \\ \\
In [IK] they prove the following \footnote{The {\itshape continuous} version of this theorem was proven by Vaserstein in [Vas]. The {\itshape algebraic} version of it (for a polynomial map $f$ of $n$ variables), apart from the trivial case $n=1$ and Cohn's well-known counterexample in case $n=k=2$ (see [C]), for $k\geq3$ and any $n$ is based upon a deep result of Suslin (see [Sus]): any matrix in $\SL_k(\C[\C^n])$ decomposes as a finite product of unipotent (and equivalently elementary) matrices.}\\ \\
{\bfseries Theorem 1.2.2.} Let $X$ be a finite dimensional reduced Stein space and $f : X \to \SL_k(\C)$ be a holomorphic mapping that is {\itshape null-homotopic}. Then there exist a natural number $N$ and holomorphic mappings $G_1,\dots, G_N : X\to \C^{k(k-1)/2}$ such that, for every $x\in X,$ $$f(x) = M_1(G_1(x))\cdots M_N(G_N(x)),$$ where, if $Z_N$ is the accordingly defined multiple variable vector, 
\normalsize
$$M_{2l}(Z_{2l})= \left(\begin{matrix} 1 & z_{12,2l} & \dots & z_{1k,2l} \\ 0 & \ddots & \ddots & \vdots  \\ \vdots & \ddots & \ddots & z_{(k-1)k,2l} \\ 0 & \dots & 0 & 1\end{matrix}\right) {\rm and} \ \ M_{2l-1}(Z_{2l-1})= \left(\begin{matrix} 1 & 0 & \dots & 0 \\ z_{21,2l-1} & \ddots & \ddots & \vdots  \\ \vdots & \ddots & \ddots & 0 \\ z_{k1,2l-1} & \dots & z_{k(k-1),2l-1} & 1\end{matrix}\right).$$
\large
So as to prove it, Ivarsson and Kutzschebauch define $\Psi_N:(\C^{k(k-1)/2})^N\to\SL_k(\C)$ as $$\Psi_N(Z_1,\dots,Z_N) = M_1(Z_1)^{-1}\cdots M_N(Z_N)^{-1}$$ and they wish to show the existence of a holomorphic map $$G = (G_1,\dots,G_N): X \to(\C^{k(k-1)/2})^N$$ such that
$$
\begin{tikzcd}[column sep=small]
 & (\C^{k(k-1)/2})^N \arrow{dr}{\Psi_N} \\
X \arrow{ur}{G} \arrow{rr}{f} && \SL_k(\C)
\end{tikzcd}
$$
is commutative.\\
Furthermore, they choose to concentrate upon 
$$
\begin{tikzcd}[column sep=small]
 & (\C^{k(k-1)/2})^N \arrow{dr}{\pi_k\circ\Psi_N} \\
X \arrow{ur}{F} \arrow{rr}{\pi_k\circ f} && \C^k\setminus\{0\}
\end{tikzcd}
$$
and to make use of the fibres of $\pi_k\circ\Psi_N=:\Phi_N,$ where $\pi_k$ denotes the projection to the $k$-th (i.e. to the last) row of a given matrix in $\SL_k(\C).$\\ \\
In this article we focus upon the case $k=2$, indicate by $n:=N\geq 3$ the number of variables involved (i.e. $n=\dim X_n+1$) and compute the homology groups of certain fibres of the maps $\Phi_n=\pi_2\circ\Psi_n$, where $\pi_2:\SL_2(\C)\to\C^2\setminus\{0\}$ is the projection to the second row. Namely, given the map $\Phi_n$, we consider, for odd $n,$ its fibre over $\{(1,a_2):\;a_2\in\C\}$ and, for even $n,$ its fibre over $\{(a_1,2):\;a_1\in\C\}$, which are biholomorphic, as an easy computation shows, to the manifolds $X_n$ we defined in section 1.1. We then compute their homology groups (see Theorem 1.1.1) and show that these fibres enjoy the algebraic volume density property (see Theorem 1.1.2). \\ \\
Another problem consists in finding the number of unipotent matrices needed in the factorization. It can be easily shown (see [IK]) that, in the {\itshape holomorphic} case as well as in the {\itshape continuous} case, there exists an upper bound depending only on the dimension $m$ of the space $X$ and the size $k$ of the matrix (while such a uniform bound does not exist in the {\itshape algebraic} case, as shown in [vdK]). As in [IK], let us indicate by $N_{\mathcal{C}}(m,k)$ the number of matrices needed to factorise any null-homotopic map from a Stein space of dimension $m$ into $\SL_k(\C)$ by continuous triangular matrices and by $N_{\mathcal{O}}(m,k)$ the number needed in the holomorphic setting. Not much about $N_{\mathcal{C}}(m,k)$ and $N_{\mathcal{O}}(m,k)$ is known so far and to determine them constitutes a very hard problem, strongly related to algebraic $K$-theory. Obviously we have $N_{\mathcal{C}} (m, k) \leq N_{\mathcal{O}}(m, k)$ and, in the case $k=2$, the proof in [IK] provides $N_{\mathcal{O}}(m,2)\leq N_{\mathcal{C}}(m,2) + 4$. The latter inequality has been considerably improved in [IKN] to $$N_{\mathcal{O}}(m,2)\leq N_{\mathcal{C,\;O}}(m,2)+2,$$ where $N_{\mathcal{C,\;O}}(m,2)$ denotes the minimal number $s$ such that every null-homotopic {\itshape holomorphic} map from a Stein space of dimension $m$ into $\SL_2(\C)$ factorises as a product of $s$ {\itshape continuous} unipotent matrices (starting with a lower triangular one). Also, in [IKN], Ivarsson and Kutzschebauch obtain some exact estimates on the number of factors, namely $$N_{\mathcal{O}}(1,2) = 4,\;\;\; N_{\mathcal{O}}(2,2) = 5.$$ Other bounds are found, in case $m=1,2$ for any $k$, by A. Brudnyi (see theorem 1.1 and proposition 1.3 of [Br] where the $K$-theory notion of {\itshape Bass stable rank} is used). We finally remark that, in order to find a bound for $N_{\mathcal{C}}(m,2),\;m\geq3,$ it would suffice to determine a topological section of the above fibration. Although such fibration is not locally trivial, but stratified locally trivial, we hope that our study of the topology of the fibres will allow us to prove the existence of a global topological section in certain cases. \\
\begin{center}
{\scshape 2. Preliminary results}
\end{center}
First of all we introduce two preliminary results which shall be used repeatedly in what follows. \\ Let us set $Y=\C^r$ and $I=(f,g)$, where $f,g\in\C^{[r]}$ are non-constant polynomials without common factor. The affine modification of $Y$ along the divisor $D=\{\overline{x}\in\C^r: f(\overline{x})=0\}\subset Y$ with centre at $C=\{f(\overline{x})=g(\overline{x})=0\}\subset D$ is the hypersurface $X\subset\C^{r+1}$ given by $X=\{(\overline{x},y)\in\C^{r+1}:f(\overline{x})y-g(\overline{x})=0\}$, obtained from $Y$ by means of the blow-up morphism $\sigma_I:X\to Y$, which is the restriction to $X$ of the natural projection $\C^{r+1}\to\C^r$, $(\overline{x},y)\to\overline{x}$. To $D$ we associate the divisor $A\subset X$ given by $A=C\times\C$.\\The following holds (see Prop. 3.1 of [KZ]): \\ \\
{\bfseries Result 1.} In the above notation, let the hypersurface $X\subset\C^{r+1}$ be the affine modification, via the projection $\sigma:X\to Y$, of $Y=\C^r$ along the divisor $D\subset Y$ with centre at $C\subset D$, to which the divisor $A=C\times\C$ is associated. Let the divisors $D$ and $A$ admit finite decompositions into irreducible components $D=\cup_{i=1}^nD_i$ and $A=\cup_{j=1}^{n'}A_j$ respectively. Suppose that $\sigma^*(D_i)=\sum_{j=1}^{n'}m_{ij}A_j$ and that $\sigma(A_j)\cap\reg D_i\neq\emptyset$ as soon as $m_{ij}>0$. Then, if the lattice vectors $b_j=(m_{1j},\dots,m_{nj})\in\Z^n,\;j=1,\dots,n',$ generate the lattice $\Z^n$, the induced map $\sigma_*:H_1(X,\Z)\to H_1(Y,\Z)=0$ is an isomorphism. \\
Assume further that there is a disjoint partition $\{1,\dots,n'\}=J_1\sqcup\dots\sqcup J_n$ such that $\sigma^*(D_i)=\sum_{j\in J_i}m_{ij}A_j\neq0,\;i=1,\dots,n.$ Let us set $d_i=\gcd(m_{ij}:j\in J_i),\;i=1,\dots,n.$ Then, if $d_1=\dots=d_n=1,$ the induced map $\sigma_*:\pi_1(X)\to\pi_1(Y)=0$ is an isomorphism.\\ \\
Also, we shall require the following result, for which we refer to [D]:\\ \\
{\bfseries Result 2.} Let $X\subset\C^{r}$ be a {\itshape smooth} hypersurface. Then $H_k(\C^r\setminus X)\cong H_{k-1}(X)$ for $k\geq1$.\\ \\
Note that, since $X$ and $\C^r\setminus X$ are connected, we have $H_0(X)\cong H_0(\C^r\setminus X)=\Z.$ \\ 
\begin{center}
{\scshape 3. Proofs}
\end{center}
\begin{center}
3.1. {\bfseries Base cases}
\end{center}
Let us prove some propositions which constitute the base cases of our induction (see section 3.2). \\ \\
{\bfseries Proposition 3.1.1.} Let $X_3=\{z_1+z_3+z_1z_3z_2=1\}\subset\C^3_{z_1,z_2,z_3}$ be the affine modification, via the projection $\sigma:X_3\to\C^2_{z_1,z_3}$, of $\C^2$ along the divisor $\Delta=\{z_1z_3=0\}\subset\C^2$ with centre at $C_3=\{z_1z_3=1-z_1-z_3=0\}=\{(1,0),(0,1)\}\subset\reg\Delta\subset\Delta$. We have $$H_0(X_3)=\Z,\;\;H_1(X_3)=0,\;\;H_2(X_3)=\Z ,\;\;H_k(X_3)=0\mbox{ for } k\geq 3.$$
{\itshape Proof.} Since $X_3$ is connected, we have $H_0(X_3)=\Z$. Note that $\C^2\setminus\Delta\simeq\C^*\times\C^*$, which is a torus.
Consider the divisors $\Delta$ and $A_3=\sigma^*(\Delta)=\sigma^{-1}(C_3)=C_3\times\C\simeq\C\sqcup\C$. Thanks to Prop. 3.1 of [KZ] one easily shows that $H_1(X_3)\cong H_1(\C^2)=0$.
As $H_*(X_3\setminus A_3)\stackrel{\sigma_*}{\cong}H_*(\C^2\setminus\Delta)$, denoting by $e$ the Euler characteristic, we have $e(X_3)=1+e(C_3)-e(\Delta)$ and, since $e(C_3)=2$ and $e(\Delta)=1$, it follows that $e(X_3)=2$ and consequently that $H_2(X_3)=\Z$, as we can have no torsion in the maximal-dimension homology group. $\;\;\diamond$\\ \\ 
{\bfseries Remark.} Consider, as in [KKAV] (par. 7), the variety $X_{p,q}=\{p(x)+q(y)+xyz=1\}\subset\C^3_{x,y,z}$, where $p$ and $q$ are polynomials such that $p(0)=q(0)=0$, $1-p(x)$ and $1-q(y)$ have simple roots only, namely $x_1,\dots,x_k$ and $y_1,\dots,y_l$. 
Thanks to Prop. 3.1 of [KZ] one easily shows that $H_1(X_{p,q})\cong H_1(\C^2)=0.$
Since we have $e(X_{p,q})=k+l$, we get $H_2(X_{p,q})=\Z^{k+l-1}$, as we can have no torsion in the maximal-dimension homology group. In particular we obtain Prop. 3.1.1 if we set $p(x)=x$, $q(y)=y$.\\ \\
{\bfseries Proposition 3.1.2.} Let $X_3^0=\{z_1+z_3+z_1z_3z_2=0\}\subset\C^3_{z_1,z_2,z_3}$ be the affine modification, via the projection $\sigma:X_3^0\to\C^2_{z_1,z_3}$, of $\C^2$ along the divisor $\Delta=\{z_1z_3=0\}\subset\C^2$ with centre at $C_3^0=\{z_1z_3=-z_1-z_3=0\}=\{(0,0)\}=\sing\Delta\subset\Delta$. We have $$H_0(X_3^0)=\Z ,\;\; H_1(X_3^0)=\Z ,\;\; H_2(X_3^0)=\Z ,\;\;H_k(X_3^0)=0\mbox{ for } k\geq 3.$$
{\itshape Proof.} Since $X_3^0$ is connected, we have $H_0(X_3^0)=\Z$. Note that $\C^2\setminus\Delta\simeq\C^*\times\C^*$, which is a torus. Consider the divisors $\Delta$ and $A_3^0=\sigma^*(\Delta)=\sigma^{-1}(C_3^0)=C_3^0\times\C=\{(0,0)\}\times\C$. Proceeding as in the proof of Prop. 3.1 of [KZ], one concludes that $H_1(X_3^0)=\Z$. As $H_*(X_3^0\setminus A_3^0)\stackrel{\sigma_*}{\cong}H_*(\C^2\setminus\Delta)$, denoting by $e$ the Euler characteristic, we have $e(X_3^0)=1+e(C_3^0)-e(\Delta)$ and, since $e(C_3^0)=e(\Delta)=1$, it follows that $e(X_3^0)=1$ and that $H_2(X_3^0)=\Z$, as we can have no torsion in the maximal-dimension homology group. $\;\;\diamond$ \\ \\
{\bfseries Proposition 3.1.3.} Let $X_4=\{z_1z_2-1+z_4(z_1+z_3+z_1z_3z_2)=0\}\subset\C^4_{z_1,z_2,z_3,z_4}$ be the affine modification, via the projection $\sigma:X_4\to\C^3_{z_1,z_2,z_3}$, of $\C^3$ along the smooth divisor $X_3^0=\{z_1+z_3+z_1z_3z_2=0\}\subset\C^3$ with centre at $C_4=\{z_1z_2-1=z_1+z_3+z_1z_2z_3=0\}=\{(z_1,1/z_1,-z_1/2):z_1\in\C^*_{z_1}\} \simeq\C^*,\; C_4\subset X_3^0$. We have $$ H_0(X_4)=\Z,\;\;H_1(X_4)\cong H_2(X_4)=0,\;\;H_3(X_4)=\Z,\;\;H_k(X_4)=0\mbox{ for } k\geq4.$$
{\itshape Proof.} Since $X_4$ is connected, we have $H_0(X_4)=\Z$. Let us consider the two long exact sequences associated to the couples $(\C^3,\C^3\setminus X_3^0)$ and $(X_4,X_4\setminus A_4)$, where $A_4=\sigma^*(X_3^0)=\sigma^{-1}(C_4)=C_4\times\C\simeq\C^*\times\C$.\\ As both divisors $X_3^0$ and $A_4$ are smooth and irreducible, thanks to Prop. 3.1 of [KZ] we can conclude that $H_1(X_4)\stackrel{\sigma_*}{\cong}H_1(\C^3)=0$.\\ Also, thanks to Thom Isomorphism theorem, 
we obtain $H_2(X_4, X_4\setminus A_4)\cong H_0(A_4)=\Z$ and $H_2(\C^3,\C^3\setminus X_3^0)\cong H_0(X_3^0)=\Z$, as well as $H_3(X_4,X_4\setminus A_4)\cong H_1(A_4)=\Z$ and $H_3(\C^3,\C^3\setminus X_3^0)\cong H_1(X_3^0)=\Z$. Thanks to the isomorphisms induced by $i\circ\sigma$, we conclude $H_2(X_4)=0$ via the Five-Lemma. Now, denoting by $e$ the Euler characteristic, we have $e(X_4)=1+e(C_4)-e(X_3^0)$. Since $e(C_4)=e(\C^*)=0$ and $e(X_3^0)=1$, it follows that $e(X_4)=0$ and consequently that $H_3(X_4)=\Z$, as we can have no torsion in the maximal-dimension homology group. $\;\;\diamond$ \\ \\
{\bfseries Proposition 3.1.4.} Let $X_5=\{z_1+z_3+z_1z_3z_2-1+z_5(z_1z_2+1+z_4(z_1+z_3+z_1z_3z_2))=0\}\subset\C^5_{z_1,\dots,z_5}$ be the affine modification, via the projection $\sigma:X_5\to\C^4_{z_1,\dots, z_4}$, of $\C^4$ along the smooth divisor $X_4^0=\{z_1z_2+1+z_4(z_1+z_3+z_1z_3z_2)=0\}\simeq\C^3\setminus X_3^0$, $X_4^0\subset\C^4,$ with centre at $C_5=\{z_1+z_3+z_1z_3z_2-1=z_1z_2+1+z_4(z_1+z_3+z_1z_3z_2)=0\}=\{(z_1,z_2,z_3,z_4):(z_1,z_2,z_3)\in X_3,z_4=-z_1z_2-1\}\simeq X_3$, $C_5\subset X_4^0.$ We have 
\normalsize
$$H_0(X_5)=\Z,\;H_1(X_5)=0,\;H_2(X_5)=\Z,\;H_3(X_5)=0,\;H_4(X_5)=\Z,\;H_k(X_5)=0\mbox{ for }k\geq 5.$$
\large
{\itshape Proof.} Here the arrow $\rightarrowtail$ will denote an {\itshape injective} map. Since $X_5$ is connected, we have $H_0(X_5)=\Z$. Let us consider the two long exact sequences associated to the couples $(\C^4,\C^4\setminus X_4^0)$ and $(X_5, X_5\setminus A_5)$, where $A_5=\sigma^*(X_4^0)=\sigma^{-1}(C_5)=C_5\times\C\simeq X_3\times\C$.\\ 
As both divisors $X_4^0$ and $A_5$ are smooth and irreducible, thanks to Prop. 3.1 of [KZ] we can conclude that $H_1(X_5)\stackrel{\sigma_*}{\cong}H_1(\C^4)=0$. By Result 2 we have $H_j(X_4^0)=\Z$ for $j=0,\dots,3$, as $X_4^0$ is connected, and also $H_j(\C^4\setminus X_4^0)=\Z$ for $j=0,\dots,3$, as $\C^4\setminus X_4^0$ is connected.
The above implies that, by Thom Isomorphism theorem, $H_j(\C^4,\C^4\setminus X_4^0)$ and $H_{j-1}(\C^4\setminus X_4^0)$ are isomorphic to $H_{j-2}(X_4^0)$ for $j=2,\dots,4$ and hence to $\Z$. \\
Thanks to the isomorphism induced by $i\circ\sigma$ we obtain $$ H_4(X_5,X_5\setminus A_5)\stackrel{\cong}{\to}H_3(X_5\setminus A_5)\stackrel{0}{\to}H_3(X_5)\stackrel{0}{\rightarrowtail}H_3(X_5,X_5\setminus A_5)=0, $$ whence $H_3(X_5)=0.$  
Observe that, by diagram commutativity, the composition $H_2(X_5)\stackrel{\tau}{\to}H_2(X_5,X_5\setminus A_5)\stackrel{\cong}{\to}H_2(\C^4,\C^4\setminus X_4^0)$, where, as above, the isomorphism is induced by $i\circ\sigma$, must be the zero map, whence $\tau=0$. It follows that $H_2(X_5\setminus A_5)\cong H_2(X_5)$, i.e. $H_2(X_5)=\Z$. \\ Denoting by $e$ the Euler characteristic, we have $e(X_5)=1+e(C_5)-e(X_4^0)=1+2-0=3$, given that $e(C_5)=e(X_3)=2$ and $ e(X_4^0)=e(\C^3\setminus X_3^0)=0$, from which it follows that $H_4(X_5)=\Z$, as we can have no torsion in the maximal-dimension homology group. $\;\;\diamond$ \\ \\
{\bfseries Proposition 3.1.5.} Let $X_6=\{z_4(z_1+z_3(z_1z_2+1))+z_1z_2-1+z_6(z_1+z_3+z_1z_3z_2+z_5(z_1z_2+1+z_4(z_1+z_3+z_1z_3z_2)))=0\}\subset\C^6_{z_1,\dots,z_6}$ be the affine modification, via the projection $\sigma:X_6\to\C^5_{z_1,\dots, z_5}$, of $\C^5$ along the smooth divisor $X_5^0=\{z_1+z_3+z_1z_3z_2+z_5(z_1z_2+1+z_4(z_1+z_3+z_1z_3z_2))=0\}\simeq\C^4\setminus X_4^0$, $X_5^0\subset\C^5,$ with centre at $C_6=\{z_4(z_1+z_3(z_1z_2+1))+z_1z_2-1=z_1+z_3+z_1z_3z_2+z_5(z_1z_2+1+z_4(z_1+z_3+z_1z_3z_2))=0\}=\{(z_1,z_2,z_3,z_4,z_5):(z_1,z_2,z_3,z_4)\in X_4,\;2z_5=-z_1-z_3-z_1z_3z_2\}\simeq X_4$, $C_6\subset X_5^0.$\\ We have
\normalsize
 $$H_0(X_6)=\Z,\;H_1(X_6)=0,\;H_2(X_6)\cong H_3(X_6)=\Z,\;H_4(X_6)=0,\;H_5(X_6)=\Z,\;H_k(X_6)=0\mbox{ for }k\geq 6.$$
\large
{\itshape Proof.} Here the arrow $\rightarrowtail$ will denote an {\itshape injective} map. Since $X_6$ is connected, we have $H_0(X_6)=\Z$. Let us consider the two long exact sequences associated to the couples $(\C^5,\C^5\setminus X_5^0)$ and $(X_6, X_6\setminus A_6)$, where $A_6=\sigma^*(X_5^0)=\sigma^{-1}(C_6)=C_6\times\C\simeq X_4\times\C$.\\ 
As both divisors $X_5^0$ and $A_6$ are smooth and irreducible, thanks to Prop. 3.1 of [KZ] we can conclude that $H_1(X_6)\stackrel{\sigma_*}{\cong}H_1(\C^5)=0$.
By Result 2 we have $H_j(X_5^0)=\Z$ for $j=0,\dots,4$, as $X_5^0$ is connected, and also $H_j(\C^5\setminus X_5^0)=\Z$ for $j=0,\dots,4$, as $\C^5\setminus X_5^0$ is connected. Hence, $H_j(X_6\setminus A_6)$ and $H_j(\C^5\setminus X_5^0)$, for sure isomorphic via $\sigma_*$, are isomorphic to $\Z$, for $j=1,\dots,3$. By Thom Isomorphism theorem, $H_3(X_6,X_6\setminus A_6)\cong H_1(A_6)=0$ and $H_4(X_6,X_6\setminus A_6)\cong H_2(A_6)=0$, so $H_3(X_6)\cong H_3(X_6\setminus A_6)$ by exactness, whence $H_3(X_6)=\Z$.\\
Thanks to the isomorphism induced by $i\circ\sigma$ we obtain $$ H_5(X_6,X_6\setminus A_6)\stackrel{\cong}{\to}H_4(X_6\setminus A_6)\stackrel{0}{\to}H_4(X_6)\stackrel{0}{\rightarrowtail}H_4(X_6,X_6\setminus A_6)=0, $$ whence $H_4(X_6)=0.$ Observe that, by diagram commutativity, the composition $H_2(X_6)\stackrel{\tau}{\to}H_2(X_6,X_6\setminus A_6)\stackrel{\cong}{\to}H_2(\C^5,\C^5\setminus X_5^0)$, where, as above, the isomorphism is induced by $i\circ\sigma$, must be the zero map, whence $\tau=0$. It follows that $H_2(X_6\setminus A_6)\cong H_2(X_6)$, i.e. $H_2(X_6)=\Z$. \\
Denoting by $e$ the Euler characteristic, we have $e(X_6)=1+e(C_6)-e(X_5^0)=1+0-1=0$, given that $e(C_6)=e(X_4)=0$ and $ e(X_5^0)=e(\C^4\setminus X_4^0)=1$, from which it follows that $H_5(X_6)=\Z$, as we can have no torsion in the maximal-dimension homology group. $\;\;\diamond$ \\ \\
\begin{center}
3.2. {\bfseries Induction}
\end{center}
This section and the following one are devoted to the proof of our theorems, which we restate along with their notation, so as to improve readability. \\ For $n\geq 3$ consider the family $\{M_n\}$ of $2\times 2$ matrices defined inductively as follows: \\

\noindent $M_3=$\;\( 
\begin{pmatrix}
1 & 0\\
z_1 & 1
\end{pmatrix}
\)
\( 
\begin{pmatrix}
1 & z_2\\
0 & 1
\end{pmatrix}
\)
\( 
\begin{pmatrix}
1 & 0\\
z_3 & 1
\end{pmatrix}
\)
=
 \( 
\begin{pmatrix}
1+z_2z_3 & z_2\\
z_1+z_3(z_1z_2+1) & z_1z_2+1
\end{pmatrix}
\),
\\ \\
$M_4=M_3$
\( 
\begin{pmatrix}
1 & z_4\\
0 & 1
\end{pmatrix}
\) =
\( 
\begin{pmatrix}
1+z_2z_3 & z_4(1+z_2z_3)+z_2\\
z_1+z_3(z_1z_2+1) & z_4(z_1+z_3(z_1z_2+1))+z_1z_2+1
\end{pmatrix}
\),
\\ \\
\dots \\ \\
$M_n=M_{n-1}$
\(
\begin{pmatrix}
1 & 0\\
z_n & 1
\end{pmatrix}
\)
\mbox{ if $n$ is odd}, \\ \\
$M_n=M_{n-1}$
\(
\begin{pmatrix}
1 & z_n\\
0 & 1
\end{pmatrix}
\)
\mbox{ if $n$ is even}.\\ \\
Also, for $n\geq 3$ consider the family $\{X_n\}$ of smooth complex algebraic hypersurfaces of $\C^n$ defined inductively as follows: 
$$X_3=\{(z_1,z_2,z_3)\in\C^3:(M_3)_{2,1}=1\}=\{p_3(z_1,z_2,z_3)=z_1+z_3+z_1z_3z_2-1=0\},$$ 
$$X_4=\{(z_1,\dots,z_4)\in\C^4:(M_4)_{2,2}=2\}=\{p_4(z_1,\dots,z_4)=z_1z_2-1+z_4(z_1+z_3+z_1z_3z_2)=0\},$$ $\dots$
$$X_n=\{(M_n)_{2,1}=1\}=\{p_n(z_1,\dots,z_n)=p_{n-2}(z_1,\dots,z_{n-2})+z_n(p_{n-1}(z_1,\dots,z_{n-1})+2)=0\}$$ if $n\geq 5$ is odd,
$$X_n=\{(M_n)_{2,2}=2\}=\{p_n(z_1,\dots,z_n)=p_{n-2}(z_1,\dots,z_{n-2})+z_n(p_{n-1}(z_1,\dots,z_{n-1})+1)=0\}$$ if $n\geq 6$ is even. Let us prove\\ \\
{\bfseries Theorem 3.2.} The following statements hold true:\\ \\
1a) $X_3\subset\C^3_{z_1,z_2,z_3}$ is the affine modification, via the projection $\sigma:X_3\to\C^2_{z_1,z_3}$, of $\C^2$ along the divisor $\Delta=\{z_1z_3=0\}\subset\C^2$ with centre at $C_3=\{z_1z_3=1-z_1-z_3=0\}=\{(1,0),(0,1)\}\subset\reg\Delta\subset\Delta$.\\ \\
1b) $X_4\subset\C^4_{z_1,z_2,z_3,z_4}$ is the affine modification, via the projection $\sigma:X_4\to\C^3_{z_1,z_2,z_3}$, of $\C^3$ along the smooth divisor $X_3^0=\{z_1+z_3+z_1z_3z_2=0\}\subset\C^3$ with centre at $C_4=\{z_1z_2-1=z_1+z_3+z_1z_2z_3=0\}=\{(z_1,1/z_1,-z_1/2):z_1\in\C^*_{z_1}\} \simeq\C^*,\; C_4\subset X_3^0$. \\ \\
1c) If $n\geq 5$ is odd $X_n\subset\C^n_{z_1,\dots,z_n}$ is the affine modification, via the projection $\sigma:X_n\to\C^{n-1}_{z_1,\dots,z_{n-1}}$, of $\C^{n-1}$ along the smooth divisor $X_{n-1}^0=\{p_{n-1}(z_1,\dots,z_{n-1})+2=0\}\subset\C^{n-1}$ with centre at $C_n=\{p_{n-2}(z_1,\dots,z_{n-2})=p_{n-1}(z_1,\dots,z_{n-1})+2=0\}\subset X_{n-1}^0$, where $X_{n-1}^0\simeq\C^{n-2}\setminus X_{n-2}^0$ and $C_n\simeq X_{n-2}$.\\ 
If $n\geq 6$ is even $X_n\subset\C^n_{z_1,\dots,z_n}$ is the affine modification, via the projection $\sigma:X_n\to\C^{n-1}_{z_1,\dots,z_{n-1}}$, of $\C^{n-1}$ along the smooth divisor $X_{n-1}^0=\{p_{n-1}(z_1,\dots,z_{n-1})+1=0\}\subset\C^{n-1}$ with centre at $C_n=\{p_{n-2}(z_1,\dots,z_{n-2})=p_{n-1}(z_1,\dots,z_{n-1})+1=0\}\subset X_{n-1}^0$, where $X_{n-1}^0\simeq\C^{n-2}\setminus X_{n-2}^0$ and $C_n\simeq X_{n-2}$.\\ \\
2a) $H_0(X_3)=\Z,\;\;H_1(X_3)=0,\;\;H_2(X_3)=\Z ,\;\;H_k(X_3)=0\mbox{ for } k\geq 3.$\\ \\
2b) $H_0(X_4)=\Z,\;\;H_1(X_4)\cong H_2(X_4)=0,\;\;H_3(X_4)=\Z,\;\;H_k(X_4)=0\mbox{ for } k\geq4.$ \\ \\
2c) If $n\geq 5$ is odd, for $0\leq j\leq n-1$ the homology groups are alternately $\Z$ or $0$, starting from $H_0(X_n)=\Z$. In particular $H_{n-2}(X_n)=0$. \\ 
If $n\geq 6$ is even, for $0\leq j\leq \frac{n-2}{2}$ the homology groups are alternately $\Z$ or $0$ starting from $H_0(X_n)=\Z$, for $\frac{n}{2}\leq j\leq n-1$ the homology groups are alternately $\Z$ or $0$ starting from $H_{n-1}(X_n)=\Z$. In particular $H_{n-2}(X_n)=0$. \\ \\
{\itshape Proof.} Statements 1a), 1b), 2a) and 2b) have been proven in propositions 3.1.1 and 3.1.3. The base cases for statements 1c) and 2c) are dealt with in propositions 3.1.4 and 3.1.5. Let us prove statement 1c) by induction on $n$. Suppose that for $n\in\N$ the statement is true and that $n+1$ is odd (so $n$ is even). We have $X_{n+1}=\{p_{n-1}+z_{n+1}(p_n+2)=0\}$ and $X_n^0=\{p_n+2=0\}$. We wish to show that $X_{n}^0\simeq\C^{n-1}\setminus X_{n-1}^0$. We may write $p_n=p_{n-2}+z_n(p_{n-1}+1)$, i.e. $X_n^0=\{p_{n-2}+2+z_n(p_{n-1}+1)=0\}$ and we are done if we show that $S=\{p_{n-2}+2=p_{n-1}+1=0\}=\emptyset$. By inductive hypothesis we have $X_{n-1}^0=\{p_{n-1}+1=0\}=\{p_{n-3}+1+z_{n-1}(p_{n-2}+2)=0\}\simeq\C^{n-2}\setminus X_{n-2}^0$, i.e. $S'=\{p_{n-3}+1=p_{n-2}+2=0\}=\emptyset$. Since we can write $p_{n-1}+1=p_{n-3}+1+z_{n-1}(p_{n-2}+2)$ we get $S=S'=\emptyset$. As regards the centre, we need to show that $C_{n+1}=\{p_{n-1}=p_n+2=0\}\simeq X_{n-1}$. We may write $p_n+2=p_{n-2}+2+z_n(p_{n-1}+1)=p_{n-2}+2+z_n$ and $C_{n+1}=\{p_{n-1}=0,\;z_n=-p_{n-2}-2\}\simeq X_{n-1}$. One proceeds similarly if $n+1$ is even. \\ 
Let us prove 2c) by induction on $n$. In what follows the arrow $\rightarrowtail$ will denote an {\itshape injective} map. Suppose that the homology of $X_{n-2}$ and the homology of $X_{n-1}$ are as wished. Let us compute the homology of $X_n$. Suppose $n$ is {\itshape odd}. \\
Since $X_n$ is connected, we have $H_0(X_n)=\Z$. Let us consider the two long exact sequences associated to the couples $(\C^{n-1},\C^{n-1}\setminus X_{n-1}^0)$ and $(X_n, X_n\setminus A_n)$, where $A_n=\sigma^{-1}(C_n)=C_n\times\C\simeq X_{n-2}\times\C$, and $X_{n-1}^0\simeq\C^{n-2}\setminus X_{n-2}^0$.\\ 
The first homology group is preserved under $\sigma$, as both divisors $X_{n-1}^0$ and $A_n$ are smooth and irreducible. Thus, thanks to Prop. 3.1 of [KZ], we can conclude that $H_1(X_n)\cong H_1(\C^{n-1})=0.$ By the same proposition we also have $\pi_1(X_n)\cong\pi_1(\C^{n-1})=0.$ One can easily show, by Result 2, that $H_k(X_{n-1}^0)=\Z\;\;\forall k\in\{0,\dots,n-2\}$ and $H_k(X_{n-1}^0)=0\;\;\forall k\geq n-1$. 
Let us consider the subsequence 
\begin{multline*}
H_{j+2}(X_n,X_n\setminus A_n)\to H_{j+1}(X_n\setminus A_n)\to H_{j+1}(X_n)\to H_{j+1}(X_n,X_n\setminus A_n)\to\\
\to H_j(X_n\setminus A_n)\to H_j(X_n)\to H_j(X_n, X_n\setminus A_n)\to H_{j-1}(X_n\setminus A_n)\to\\ 
\to H_{j-1}(X_n)\to H_{j-1}(X_n, X_n\setminus A_n)\to H_{j-2}(X_n\setminus A_n)\to H_{j-2}(X_n)
\end{multline*}
which may be written as
\begin{multline*}
H_j(X_{n-2})\to\Z\to H_{j+1}(X_n)\to H_{j-1}(X_{n-2})\to\\
\to\Z\to H_j(X_n)\to H_{j-2}(X_{n-2})\to\Z\to\\ 
\to H_{j-1}(X_n)\to H_{j-3}(X_{n-2})\to\Z\to H_{j-2}(X_n).
\end{multline*}
{\bfseries Case 1.} Suppose $H_{j-3}(X_{n-2})=0.$ Then, by inductive hypothesis, $H_{j-2}(X_{n-2})=\Z$, $H_{j-1}(X_{n-2})=0$ and $H_j(X_{n-2})=\Z$. As $H_{j-2}(X_{n-2})\cong H_{j-2}(X_{n-1}^0)$ via $(i\circ\sigma)_*$ the subsequence becomes 
\begin{multline*}
\Z\to\Z\to H_{j+1}(X_n)\to 0
\to\Z\to H_j(X_n)\to \Z\stackrel{\cong}{\to}\Z 
\stackrel{0}{\to} H_{j-1}(X_n)\stackrel{0}{\rightarrowtail} 0\to\Z\to H_{j-2}(X_n),
\end{multline*}
whence $H_{j-1}(X_n)=0.$
Observe that, by diagram commutativity, the composition $H_j(X_n)\stackrel{\tau}{\to}\Z\stackrel{\cong}{\to}H_{j-2}(X_{n-1}^0)$, where, as above, the isomorphism is induced by $i\circ\sigma$, must be the zero map, whence $\tau=0$. It follows that $H_j(X_n\setminus A_n)\cong H_j(X_n)$, i.e. $H_j(X_n)=\Z$. Proceeding as we did for $H_{j-1}(X_n)$, we get $H_{j+1}(X_n)=0.$\\
{\bfseries Case 2.} Suppose $H_{j-3}(X_{n-2})=\Z.$ Then, by inductive hypothesis, $H_{j-2}(X_{n-2})=0$, $H_{j-1}(X_{n-2})=\Z$ and $H_j(X_{n-2})=0.$ Proceeding as we did for $H_j(X_n)$ in case 1, we get $H_{j-1}(X_n)=\Z$. Proceeding as we did for $H_{j-1}(X_n)$ in case 1, we get $H_{j}(X_n)=0.$ Again, proceeding as we did for $H_j(X_n)$ in case 1, we get $H_{j+1}(X_n)=\Z.$\\ The above shows the alternation of the homology groups, as wished. \\
We remark that an alternative way of obtaining $H_{n-1}(X_n)=\Z$ when $n$ is odd consists in noticing that, as an easy induction shows, $e(X_n)=1+e(X_{n-2})-e(X_{n-1}^0)=1+(n-2-\frac{n-3}{2})-0=\frac{n+1}{2}$, where $e$ denotes the Euler characteristic. This ends the proof of point 2c) in the case in which $n$ is {\itshape odd}. \\
Let us suppose $n$ is {\itshape even}. Since $X_n$ is connected, we have $H_0(X_n)=\Z$. Let us consider the two long exact sequences associated to the couples $(\C^{n-1},\C^{n-1}\setminus X_{n-1}^0)$ and $(X_n, X_n\setminus A_n)$, where $A_n=\sigma^{-1}(C_n)=C_n\times\C\simeq X_{n-2}\times\C$, and $X_{n-1}^0\simeq\C^{n-2}\setminus X_{n-2}^0$.\\ 
The first homology group is preserved under $\sigma$, as both divisors $X_{n-1}^0$ and $A_n$ are smooth and irreducible. Thus, thanks to Prop. 3.1 of [KZ], we can conclude that $H_1(X_n)\cong H_1(\C^{n-1})=0.$ By the same proposition we also have $\pi_1(X_n)\cong\pi_1(\C^{n-1})=0.$ One can easily show, by Result 2, that $H_k(X_{n-1}^0)=\Z\;\;\forall k\in\{0,\dots,n-2\}$ and $H_k(X_{n-1}^0)=0\;\;\forall k\geq n-1$. The long exact sequence of $(X_n,X_n\setminus A_n)$ can be subdivided as follows into three subsequences:
\begin{multline}
H_n(X_n,X_n\setminus A_n)\to\dots\to H_{1+\frac{n}{2}}(X_n,X_n\setminus A_n),
\end{multline}
\begin{multline}
H_{1+\frac{n}{2}}(X_n,X_n\setminus A_n)\to H_{\frac{n}{2}}(X_n\setminus A_n)\to H_{\frac{n}{2}}(X_n)\to H_{\frac{n}{2}}(X_n,X_n\setminus A_n)\to\\
\to H_{\frac{n-2}{2}}(X_n\setminus A_n)\to H_{\frac{n-2}{2}}(X_n)\to H_{\frac{n-2}{2}}(X_n,X_n\setminus A_n),
\end{multline}
\begin{multline}
H_{\frac{n-2}{2}}(X_n,X_n\setminus A_n)\to\dots\to H_2(X_n,X_n\setminus A_n).
\end{multline}
In sequence $(1)$ the block $0\to\Z\to H_j(X_n)\to\Z$ alternates with the block $\Z\to\Z\to H_{j-1}(X_n)\to 0$, where the middle group is in common. Each of these blocks can be solved on its own. Namely, by diagram commutativity, the composition $H_j(X_n)\stackrel{\tau}{\to}\Z\stackrel{\cong}{\to}H_{j-2}(X_{n-1}^0)$, where the isomorphism is induced by $i\circ\sigma$, must be the zero map, whence $\tau=0$. It follows that $H_j(X_n\setminus A_n)\cong H_j(X_n)$, i.e. $H_j(X_n)=\Z$. Also, as $H_{j-2}(X_{n-2})\cong H_{j-2}(X_{n-1}^0)$ via $(i\circ\sigma)_*$ the second block becomes $$
\Z\stackrel{\cong}{\to}\Z\stackrel{0}{\to}H_{j-1}(X_n)\stackrel{0}{\rightarrowtail} 0.$$ whence $H_{j-1}(X_n)=0.$\\
In sequence $(3)$ the block $\Z\to\Z\to H_l(X_n)\to 0$ alternates with the block $0\to\Z\to H_{l-1}(X_n)\to\Z$, where the middle group is in common. Each of these blocks can be solved on its own like the ones above. Thus we have solved the upper and the lower part of the long exact sequence associated to $(X_n,X_n\setminus A_n)$. \\
So as to solve subsequence $(2)$ we consider the following two cases.\\ 
{\bfseries Case 1.} Suppose that $H_{\frac{n-4}{2}}(X_{n-2})\cong H_{\frac{n-2}{2}}(X_{n-2})=0.$ Subsequence $(2)$ becomes $$0=H_{\frac{n-2}{2}}(X_{n-2})\to\Z\stackrel{\cong}{\to}H_{\frac{n}{2}}(X_n)\to H_{\frac{n-4}{2}}(X_{n-2})=0\to\Z\to H_{\frac{n-2}{2}}(X_n)\to H_{\frac{n-4}{2}-1}(X_{n-2}),$$ whence $H_{\frac{n}{2}}(X_n)=\Z$. As, by inductive hypothesis, $H_{\frac{n-4}{2}-1}(X_{n-2})=\Z$, proceeding as we did for the first block of sequence $(1)$, we obtain $H_{\frac{n-2}{2}}(X_n)=\Z$. Hence\\ $H_{\frac{n}{2}}(X_n)\cong H_{\frac{n-2}{2}}(X_n)=\Z.$ \\
{\bfseries Case 2.} Suppose that $H_{\frac{n-4}{2}}(X_{n-2})\cong H_{\frac{n-2}{2}}(X_{n-2})=\Z.$ Subsequence $(2)$ becomes $$\Z=H_{\frac{n-2}{2}}(X_{n-2})\to\Z\to H_{\frac{n}{2}}(X_n)\to H_{\frac{n-4}{2}}(X_{n-2})=\Z\to\Z\to H_{\frac{n-2}{2}}(X_n)\to H_{\frac{n-4}{2}-1}(X_{n-2}).$$ Thanks to the isomorphisms induced by $i\circ\sigma$, we conclude $H_{\frac{n}{2}}(X_n)=0$ via the Five-Lemma. As, by inductive hypothesis, $H_{\frac{n-4}{2}-1}(X_{n-2})=0$, proceeding as we did for the second block of sequence $(1)$, we obtain $H_{\frac{n-2}{2}}(X_n)=0$. Hence $H_{\frac{n}{2}}(X_n)\cong H_{\frac{n-2}{2}}(X_n)=0.$ We remark that an alternative way of obtaining $H_{n-1}(X_n)=\Z$ when $n$ is even consists in noticing that, as an easy induction shows, $e(X_n)=1+e(X_{n-2})-e(X_{n-1}^0)=1+0-1=0$, where $e$ denotes the Euler characteristic. This concludes the proof. $\;\;\diamond$ \\ \\
\begin{center}
3.3. {\bfseries Volume forms}
\end{center}
We prove the following \\ \\
{\bfseries Lemma 3.3.1.} For each $n\geq 3$ the variety $X_n=\{p_n=0\}$ can be equipped with an algebraic volume form $\omega$.\\ \\
{\itshape Proof.} Let us express each variable in terms of the others. We obtain $z_i=\zeta_i(z_1,\dots,\widehat{z_i},\dots,z_n)$. Let us consider the cover $\{U_i\}_{i=1}^n$ of $X_n$ with $U_i=X_n\setminus K_i$, where $K_i=\{\underline{z}\in X_n:\partial_{z_i}p_n(\underline{z})=\frac{\partial p_n}{\partial z_i}(\underline{z})=0\}$, and the set $\{(\omega_i,U_i)\}_{i=1}^n$, where $\omega_i=\frac{1}{\partial_{z_i}p_n}\;dz_1\wedge\dots\wedge\widehat{dz_i}\wedge\dots\wedge dz_n$ is a volume form defined on $\varphi_i^{-1}(U_i)=\C^{n-1}_{z_1,\dots,\widehat{z_i},\dots,z_n}\setminus H_i$, $H_i=\{\underline{z}\in\C^{n-1}:\partial_{z_i}p_n(\underline{z})=0\}\subset\C^{n-1}_{z_1,\dots,\widehat{z_i},\dots,z_n}$ and $\varphi_{i}:\C^{n-1}_{z_1,\dots,\widehat{z_i},\dots,z_n}\setminus H_i\to U_i$ is defined via $\varphi_i(z_1,\dots,\widehat{z_i},\dots,z_n)=(z_1,\dots,z_i=\zeta_i(z_1,\dots,\widehat{z_i},\dots,z_n),\dots,z_n)$, while $\varphi_i^{-1}:U_i\to \C^{n-1}_{z_1,\dots,\widehat{z_i},\dots,z_n}\setminus H_i$ takes $(z_1,\dots,z_n)$ to $(z_1,\dots,\widehat{z_i},\dots,z_n)$. For $i\neq j$ the map $\varphi_i^{-1}\circ\varphi_j:\varphi_j^{-1}(U_i\cap U_j)\to\varphi_i^{-1}(U_i\cap U_j)$, is such that, up to sign, $(\varphi_i^{-1}\circ\varphi_j)^*\omega_i=\omega_j$. Hence the local forms $(\omega_i,U_i)$ are pairwise compatible and, glued together, constitute a volume form $\omega$ on $X_n$. $\;\;\diamond$ \\ \\
Let $\omega$ be an algebraic volume form on $X_n$. Let us consider, for $i\neq j\in\{1,\dots,n\}$, the following complete $\omega$-divergence-free algebraic vector fields $\delta_{ij}$ on $X_n$: $$\delta_{ij}=\frac{\partial p_n}{\partial z_i}\;\frac{\partial}{\partial z_j}-\frac{\partial p_n}{\partial z_j}\;\frac{\partial}{\partial z_i}.$$ 
So as to prove that $X_n$ has the algebraic volume density property, it suffices to show, since for each $n\geq 3$ we have $H_{n-2}(X_n)=0$, that $$\Theta(\Lie_{\alg}^{\omega}(X_n))\supset\mathcal{B}_{n-2}(X_n),$$ where $\Theta:\AVF_{\omega}(X_n)\to\mathcal{Z}_{n-2}(X_n)$ is the isomorphism sending the $\omega$-divergence-free algebraic vector field $\xi$ to the closed $(n-2)$-form $\iota_{\xi}\omega$, defined by means of the interior product $\iota$ of $\xi$ and $\omega$, $\Lie_{\alg}^{\omega}(X_n)$ is the Lie algebra generated by the set $\IVF_{\omega}(X_n)$ of complete $\omega$-divergence-free algebraic vector fields on $X_n$ and $\mathcal{B}_{n-2}(X_n)$ is the space of exact algebraic $(n-2)$-forms on $X_n$. \\ Now, the above inclusion does hold, as an easy proof shows.  Let us go through its key steps for arbitrary $n\geq 3.$ \\ 
First of all, by Grothendieck's theorem we have $\mathcal{B}_{n-2}(X_n)\subset\Theta(\AVF_{\omega}(X_n))$, i.e. for every algebraic $(n-3)$-form $\alpha$ on $X_n$ there exists $\xi\in\AVF_{\omega}(X_n)$ such that $\Theta(\xi)=\iota_{\xi}\omega={\rm d}\alpha$. It suffices to show $\xi$ can be chosen in $\Lie_{\alg}^{\omega}(X_n)$. Each coefficient of $\alpha$ is a regular function on $X_n$, thus being the sum of monomials of the kind $z_{i_1}^{r_1},\; z_{i_1}^{r_1}z_{i_2}^{r_2},\dots,\;z_{i_1}^{r_1}z_{i_2}^{r_2}\cdots z_{i_{n-1}}^{r_{n-1}}$, where $i_1<i_2<\;\cdots\;<i_{n-1}$ are elements of the set $\{1,2,\dots,n\}$. We further have that the isomorphism $\Lambda$ induced by the map $\alpha\to\xi$ is such that, up to constant factors, if $i,j\neq1$ (otherwise we obtain a similar result), $\Lambda^{-1}(\delta_{ij})=z_1\;dz_2\;\wedge\;\cdots\;\wedge\;\widehat{dz_i}\;\wedge\;\cdots\;\wedge\;\widehat{dz_j}\;\wedge\;\cdots\;\wedge\; dz_n$, as $\Ker(\delta_{ij})=\C[z_1,\dots,\widehat{z_{i}},\dots,\widehat{z_j},\dots,z_n]$. The other summands of the kind $z_{i_1}$ and $z_{i_1}^{r_1}$ are dealt with similarly. The Lie bracket $[\delta_{ij},\delta_{kl}]\in\Lie_{\alg}^{\omega}(X_n),\; i\neq j,k\neq l$ is sent to a $(n-3)$-form whose coefficients are monomials $z_{i_1}z_{i_2}$. Monomials $z_{i_1}^{r_1}z_{i_2}^{r_2}$ are treated similarly. As regards summands of the kind $z_{i_1}z_{i_2}z_{i_3},$ without loss of generality we can suppose that the extra variable is not in the kernel of both vector fields $\delta_{ij}$ and $\delta_{kl}$. Then the product of the extra variable with the vector field annihilating it, bracketed with the other vector field belongs to $\Lie_{\alg}^{\omega}(X_n)$ and is sent to a $(n-3)$-form whose coefficients are monomials in three variables, linear in each variable. One proceeds similarly for $z_{i_1}z_{i_2}\cdots z_{i_m}$ and $z_{i_1}^{r_1}z_{i_2}^{r_2}\cdots z_{i_m}^{r_m},\; 4\leq m\leq n-1$. Hence $\xi\in\Lie_{\alg}^{\omega}(X_n)$ and the following holds true \\ \\
{\bfseries Theorem 3.3.2.} For each $n\geq 3$ the variety $X_n\subset\C^n$ is such that $\Theta(\Lie_{\alg}^{\omega}(X_n))\supset\mathcal{B}_{n-2}(X_n)$, i.e. for every algebraic $(n-3)$-form $\alpha$ on $X_n$ there exists $\xi\in\Lie_{\alg}^{\omega}(X_n)$ such that $\Theta(\xi)=\iota_{\xi}\omega={\rm d}\alpha$. Hence $\Lie_{\alg}^{\omega}(X_n)=\AVF_{\omega}(X_n)$, i.e. $X_n$ has the algebraic volume density property.

\newpage
\small
\noindent
{\bfseries Bibliography}\\ \\
1. [A]:  E. Anders\'{e}n, {\itshape Volume Preserving Automorphisms of $\C^n$}, Complex Variables 14, 223-235 (1990); \\ \\
2. [AL]: E. Anders\'{e}n, L. Lempert, {\itshape On The Group of Holomorphic Automorphisms of $\C^n$}, Invent. Math. 110, 371-388 (1992); \\ \\
3. [AF]: R. Andrist, F. Forstneri\u{c}, T. Ritter, E. F. Wold, {\itshape Proper holomorphic embeddings into Stein manifolds with the density property}, J. Anal. Math., 130:135-150, 2016;\\ \\
4. [AW]: R. B. Andrist, E. F. Wold, {\itshape Riemann surfaces in Stein manifolds with the density property}, Ann. Inst. Fourier (Grenoble), 64(2):681-697, 2014;\\ \\
5. [Br]: A. Brudnyi, {\itshape On the Bass stable rank of Stein Algebras}, arXiv:1803.10790v1 [math.CV] 28 Mar. 2018; \\ \\ 
6. [C]: P. M. Cohn, {\itshape On the structure of the $\GL_2$ of a ring}, Inst. Hautes \'{E}tudes Sci. Publ. Math. 30 (1966), 5-53. MR 0207856. Zbl 0144.26301. http: //dx.doi.org/10.1007/BF02684355; \\ \\
7. [D]: A. Dimca, A. N\'{e}methi, {\itshape Hypersurface complements, Alexander modules and monodromy}, arXiv:math/0201291v1 [math.AG] 29 Jan. 2002;\\ \\
8. [F]: F. Forstneri\u{c}, {\itshape Holomorphic embeddings and immersions of Stein manifolds: a survey}, arXiv:1709.05630v5  [math.CV]  21 Jul. 2018; \\ \\
9. [F2]: F. Forstneri\u{c}, {\itshape Proper holomorphic immersions into Stein manifolds with the density property}, ArXiv e-prints, Mar. 2017. https://arxiv.org/abs/1703.08594;\\ \\
10. [FH]: F. Forstneri\u{c}, {\itshape Stein manifolds and holomorphic mappings. The homotopy principle in complex analysis} (2nd edn), volume 56 of {\itshape Ergebnisse der Mathematik und ihrer Grenzgebiete}, 3. Folge, Berlin, Springer, 2017; \\ \\
11. [G]: M. Gromov, {\itshape Oka's principle for holomorphic sections of elliptic bundles}, J. Amer. Math. Soc. 2 (1989), 851-897. MR 1001851. Zbl 0686.32012. http://dx.doi.org/10.2307/1990897;\\ \\
12. [IK]: B. Ivarsson, F. Kutzschebauch, {\itshape Holomorphic factorization of mappings into $\SL_n(\C)$}, Annals of Mathematics 175 (2012), 45--69; \\ \\
13. [IKN]: B. Ivarsson, F. Kutzschebauch, {\itshape On the number of factors in the unipotent factorization of holomorphic mappings into $\SL_2(\C)$}, Proc. Amer. Math. Soc., on June 23, 2011, http://dx.doi.org/10.1090/S0002-9939-2011-11025-6; \\ \\ 
14. [K]: F. Kutzschebauch, {\itshape Manifolds with infinite dimensional group of holomorphic automorphisms and the linearization problem}, arXiv:1903.00970v1 [math.CV] 3 Mar. 2019; \\ \\
15. [KKAV]: Sh. Kaliman, F. Kutzschebauch, {\itshape On algebraic volume density property}, Springer Science+Business Media, New York (2015);\\ \\
16. [KKAVDP]: Sh. Kaliman, F. Kutzschebauch, {\itshape Algebraic volume density property of affine algebraic manifolds}, Invent. math. (2010) 181: 605. https://doi.org/10.1007/s00222-010-0255-x; \\ \\
17. [KKH]: Sh. Kaliman, F. Kutzschebauch, {\itshape Algebraic (volume) density property for affine homogeneous spaces}, Math. Ann. (2017) 367: 1311. https://doi.org/10.1007/s00208-016-1451-9; \\ \\
18. [KKPS]: Sh. Kaliman, F. Kutzschebauch, {\itshape On the Present State of the Anders\'{e}n-Lempert Theory},  Centre de Recherches Math\'{e}matiques CRM, Proceedings and Lecture Notes, 2011; \\ \\
19. [KZ]: Sh. Kaliman, M. Zaidenberg, {\itshape Affine modifications and affine hypersurfaces with a very transitive automorphism group}, arXiv:math/9801076v2 [math.AG], 27 Feb. 1998;\\ \\
20. [N]: R. Narasimhan, {\itshape Imbedding of holomorphically complete complex spaces}, Amer. J. Math.82 (1960), 917-934. MR0148942 (26:6438); \\ \\
21. [V]: D. Varolin, {\itshape A general notion of shears, and applications}, Michigan Math. J. 46 (1999), no. 3, 533-553. MR 1721579; \\ \\
22. [V1]: D. Varolin, {\itshape The density property for complex manifolds and geometric structures}, J. Geom. Anal. (2001) 11: 135. https://doi.org/10.1007/BF02921959; \\ \\
23. [Vas]: L. N. Vaserstein, {\itshape Reduction of a matrix depending on parameters to a diagonal form by addition operations}, Proc. Amer. Math. Soc. 103 (1988), 741-746. MR 0947649. Zbl 0657.55005. http://dx.doi.org/10.2307/ 2046844; \\ \\
24. [vdK]: W. van der Kallen, {\itshape $\SL_3(\C[X])$ does not have bounded word length}, in Algebraic K-theory, Part I (Oberwolfach, 1980), Lecture Notes in Math. 966, Springer-Verlag, New York, 1982, pp. 357-361. MR 0689383. Zbl 0935.20501. http://dx.doi.org/10.1007/BF02684355; \\ \\
25. [P]: A. Ramos-Peon, {\itshape Non-algebraic examples of manifolds with the volume density property}, arXiv:1602.07862v1  [math.CV]  25 Feb. 2016; \\ \\
26. [R]: J.--P. Rosay, {\itshape Automorphism of $\C^n,$ a survey of Anders\'{e}n--Lempert theory and applications}, Complex Geometric Analysis in Pohang (Pohang, 1997);\\ \\
27. [Sus]: A. A. Suslin, {\itshape The structure of the special linear group over rings of polynomials}, Izv. Akad. Nauk SSSR Ser. Mat. 41 (1977), 235-252, 477. MR 0472792. Zbl 0378.13002; \\ \\
28. [T]: T. Ritter, {\itshape A strong Oka principle for embeddings of some planar domains into $\C\times\C^*$}, J. Geom. Anal. 23, p. 571-597 (2013);\\ \\
29. [T1]: T. Ritter, {\itshape Acyclic embeddings of open Riemann surfaces into new examples of elliptic manifolds}, Proc. Amer. Math. Soc. 141, p. 597-603 (2013).
\end{document}